\documentclass{amsart}
\usepackage{mathrsfs,amsthm,amscd,amssymb,latexsym,multirow,graphics,xypic}

\newcommand{\COMPACT}{\mathscr{K}}

\newcommand{\DIAGRAM}{\mathscr{D}}

\newcommand{\GENERIC}{\mathscr{C}}

\newcommand{\STREAMS}{\mathscr{S}}
\newcommand{\SPACES}{\mathscr{T}}

\newcommand{\colim}{\mathrm{colim}}

\newcommand{\graph}[1]{\mathrm{graph}(#1)}
\newcommand{\id}{\mathrm{id}}
\newcommand{\domain}[1]{\mathrm{domain}(#1)}

\newcommand{\upper}[2]{\leqslant_{#1}\!\![#2]}
\newcommand{\downer}[2]{\leqslant_{#1}^{-1}\!\![#2]}

\newcommand{\BASIS}{\mathscr{B}}

\newtheorem{thm}{Theorem}[section]
\newtheorem*{thm:embed}{Theorem 4.7}
\newtheorem*{thm:xclosed}{Theorem 5.12}

\newtheorem{lem}[thm]{Lemma}

\newtheorem{prop}[thm]{Proposition}

\newtheorem{defn}[thm]{Definition}

\newtheorem{eg}[thm]{Example}

\begin{document}                                                                                   
\title{A convenient category of locally preordered spaces} 
\author{Sanjeevi Krishnan}  
\address{Laboratoire d'Informatique de l'\'{E}cole Polytechnique\\Palaiseau, France}
\begin{abstract}
	As a practical foundation for a homotopy theory of abstract spacetime, we extend a category of certain compact partially ordered spaces to a convenient category of ``locally preordered'' spaces.
	In particular, we show that our new category is Cartesian closed and that the forgetful functor to the category of compactly generated spaces creates all limits and colimits.
\end{abstract}
\maketitle
\section{Introduction}\label{sec:intro}
A homotopy theory which respects the flow of time on a machine state space $X$ can detect behavior unseen by the classical homotopy theory of $X$, as shown in \cite{fgr:ditop, fgrh:components, gh:components2, grandis:images, pratt:models}.  
Take as an example $X=\mathbb{S}^1$, the state space of a cyclical process.  
We might write $x\leqslant_{\mathbb{S}^1}y$ to express the reachability of state $y$ from state $x$, but the resulting preorder $\leqslant_{\mathbb{S}^1}$ would have graph $\mathbb{S}^1\times\mathbb{S}^1$ and therefore would fail to distinguish between clockwise and counterclockwise traversals of $\mathbb{S}^1$.  
A single preorder need not describe the ``local'' behavior of time.

The literature adopts several distinct formalisms to encode the time in abstract spacetime: an ``atlas'' on $X$ of partial orders in \cite{fgr:ditop}, a distinguished choice of paths on $X$ in \cite{grandis:d}, a quotient map $E\rightarrow X$ of spaces together with a preorder on $E$ in \cite{grandis:ineq}, and structure maps ``almost'' turning $X$ into an internal topological category in \cite{gaucher:flows}; the resulting categories of abstract spacetime share certain characteristic features identified in \cite{haucourt:dicategories}.
We propose an alternative axiomatization of abstract spacetime which generalizes the partially ordered spaces we encounter in nature while forming a category convenient for a homotopy theorist.

As in \cite{krishnan:thesis}, we define a \textit{stream} $(X,\leqslant)$ in \S\ref{subsec:definitions} to be a space $X$ equipped with a \textit{circulation} $\leqslant$, a coherent preordering $U\mapsto\leqslant_U$ of its open subsets.
The category $\STREAMS$ of maps preserving all structure in sight is complete and cocomplete.   
Colimits, limits, and \textit{substreams} are colimits, limits, and subspaces equipped with universal circulations, which we construct as \textit{pushforwards} and \textit{cosheafifications} of \textit{pullbacks} in \S\ref{subsec:cosheafafications}.
We can think of $\STREAMS$ as an extension of the category $\COMPACT$ of locally convex, compact Hausdorff partially ordered spaces whose bounded intervals are closed and connected.

\begin{thm:embed}
	There exists a full and concrete embedding
	\begin{equation}\label{eqn:intro.embedding}
		\COMPACT\hookrightarrow\STREAMS
	\end{equation}
	sending each $\COMPACT$-object to a unique stream it underlies.  
	The image of (\ref{eqn:intro.embedding}) contains all compact Hausdorff streams having locally convex underlying preordered spaces whose bounded intervals are closed.
\end{thm:embed}

The category $\STREAMS$ is not Cartesian closed.
In \S\ref{sec:compactly.flowing}, we replace $\STREAMS$ with its full subcategory $\STREAMS'$ of \textit{compactly flowing} streams, analogous to the category $\SPACES'$ of compactly generated spaces.
A stream is compactly flowing if it is locally compact Hausdorff, for example.
Our new forgetful functor $\STREAMS'\rightarrow\SPACES'$ creates limits and colimits as before, our old embedding $\COMPACT\hookrightarrow\STREAMS$ corestricts to a new embedding $\COMPACT\hookrightarrow\STREAMS'$, the formation of products and colimits in $\STREAMS'$ remains intuitive, and $\STREAMS'$ contains all streams of interest - but now we have an additional convenience.

\begin{thm:xclosed}
	The category $\STREAMS'$ is Cartesian closed.
\end{thm:xclosed}

We translate our formalism into others where topical.  
Examples \ref{eg:d-spaces}, \ref{eg:dipaths}, \ref{eg:d.beaucoup}, and \ref{eg:k.d-spaces} compare streams with the d-spaces of \cite{grandis:d}, while Examples \ref{eg:atlases}, \ref{eg:quotients}, and \ref{eg:local.embed} compare streams with the locally partially ordered spaces of \cite{fgr:ditop}.    
The comparison functors in the examples facilitate the construction of abstract spacetime in all three settings.  
We suggest a possible line of research in \S\ref{sec:future}. 

\section{Preordered spaces}\label{sec:preordered.spaces}

We fix some order-theoretic notation in \S\ref{subsec:order} and recall the basic definitions of preordered spaces in \S\ref{subsec:preordered.spaces}.

\subsection{Some order-theoretic conventions}\label{subsec:order}

Recall that a \textit{relation $R$ on a set $X$}, generalizing a function $X\rightarrow X$, encapsulates the data of its \textit{domain} $\domain{R}=X$ and its \textit{graph} $\graph{R}$, a subset of $X\times X$.  

\begin{eg}
	A relation $f$ on a set $X$ is a function $f:X\rightarrow X$ if for each $x\in X$, there exists a unique $y\in Y$ such that $(x,y)\in\graph{f}$.		
\end{eg}

For a relation $R$, write $x\;R\;y$ if $(x,y)\in\graph{R}$.  
Write
$$x_0\;R_1\;x_1\;R_2\;\ldots\;x_{n-1}\;R_n\;x_n$$
for a sequence $R_1,\ldots,R_n$ of relations if $x_{i-1}\;R_i\;x_i$ for each $0<i\leq n$.
Certain operations on functions $X\rightarrow X$ generalize to arbitrary relations.

\begin{defn}
	Consider a relation $R$ on a set $X$.  
	For each subset $A\subset X$, define $R_{\restriction A}$ to be the relation on $A$ having graph
	$$\graph{R}\cap (A\times A).$$
	
	Define $R^{-1}$ to be the relation on $X$ having graph
	$$\{(y,x)\;|\;(x,y)\in\graph{R}\}.$$
	
	For each $x\in X$, define $R[x]$ to be the subset $\{y\;|\;x\;R\;y\}\subset X$.
\end{defn}

A major reason we choose to distinguish a relation from its graph is so that we unambiguously can denote ``product relations,'' generalizations of products of functions $X\rightarrow X$.

\begin{defn}
	Consider an indexing set $\mathcal{I}$ and a relation $R_i$ on a set $X_i$ for each $i\in\mathcal{I}$.  
	For each $j\in\mathcal{I}$, let $\pi_j:\prod_{i\in\mathcal{I}}X_i\rightarrow X_j$ denote projection.
	Define $\prod_{i\in\mathcal{I}}R_i$ to be the relation on $\prod_{i\in\mathcal{I}}X_i$ having graph
	$$\bigcap_{i\in\mathcal{I}}\{(x,y)\in (\prod_{i\in\mathcal{I}}X_i)\times(\prod_{i\in\mathcal{I}}X_i)\;|\;\pi_i(x)\;R_i\;\pi_i(y)\}.$$
	If $\mathcal{I}$ consists of two elements, say $0$ and $1$, write $R_0\times R_1$ for $\prod_{i\in\mathcal{I}}R_i$.
\end{defn}

A relation $\leqslant_X$ is a \textit{preorder} if $x\leqslant_Xx$ for all $x$ and $x\leqslant_Xz$ whenever $x\leqslant_Xy\leqslant_Xz$.  
Furthermore, a preorder $\leqslant_X$ is a \textit{partial order} if $x=y$ whenever $x\leqslant_Xy\leqslant_Xx$.
A \textit{preordered set} $(X,\leqslant_X)$ is a set $X$ equipped with a preorder $\leqslant_X$ on $X$.

\begin{eg}
	The identity function $\mathrm{id}_X:X\rightarrow X$ on a set $X$ is the preorder on $X$ with smallest graph.   
\end{eg}

\begin{eg}
	The standard order $\leqslant_{\mathbb{I}}$ on the unit interval $\mathbb{I}=[0,1]$ is a partial order.  
\end{eg}

We can generalize the closed intervals and convex subsets of $\mathbb{R}$ to the setting of arbitrary preordered sets.  

\begin{defn}
	Consider a preordered set $(X,\leqslant_X)$.
	The \textit{bounded intervals} of $(X,\leqslant_X)$ are all subsets of $X$ of the form
	$$\upper{X}{x}\;\cap\downer{X}{y},\quad\quad x,y\in X.$$
	More generally, a subset $C\subset X$ is \textit{convex in $(X,\leqslant_X)$} if $y\in C$ whenever $x\leqslant_Xy\leqslant_Xz$, for all $x,z\in C$ and all $y\in X$.
\end{defn}

\begin{eg}[Geometric convexity]
	Consider a real vector space $V$ and let $x,y,v$ denote points in $V$.
	For each $v$, define a relation $R(v)$ on $V$ by the rule $x\;R(v)\;y$ if $y-x=\lambda v$ for some non-negative scalar $\lambda$.
	We leave it to the reader to check that a subset $A\subset V$ is convex in the usual sense if and only if $A$ is convex in $(V,R(v))$ for each $v$.
\end{eg}

We denote the \textit{transitive-reflexive closure} of a relation $R$ on a set $X$, the preorder on $X$ with smallest graph containing $\graph{R}$, by the notation $R^\infty$ in the following statements.  
We omit proofs, as they are straightforward.

\begin{lem}\label{lem:transitive.closure}
	For all functions $f:X\rightarrow Y$ and preorders $\leqslant_Y$ on $Y$, 
	$$(f\times f)(\graph{R^\infty})\subset\graph{\leqslant_Y}$$
	for every relation $R$ whose graph maps into $\graph{\leqslant_Y}$ under $f\times f$.
\end{lem}

\begin{lem}
	For all relations $R$ and $S$, $(R\times S)^\infty=R^\infty\times S^\infty$. 
\end{lem}

We single out transitive-reflexive closures of ``unions.''

\begin{defn}
	Consider a family $\{(A,\leqslant_A)\}_{A\in\mathcal{O}}$ of preordered sets.
	Define $\bigvee_{A\in\mathcal{O}}\leqslant_A$ to be the transitive-reflexive closure of the relation on $\bigcup\mathcal{O}$ with graph $\bigcup_{A\in\mathcal{O}}\graph{\leqslant_A}$.
\end{defn}

The operation $\times$ distributes over $\bigvee$ in the sense of the following lemma, whose proof is straightforward.

\begin{lem}\label{lem:demorgan}
	For all sets $\{\leqslant_X\}\cup\{(\leqslant_A)\}_{A\in\mathcal{O}}$ of preorders, 
	$$\leqslant_X\times(\bigvee_{A\in\mathcal{O}}\leqslant_A)=\bigvee_{A\in\mathcal{O}}(\leqslant_X\times\leqslant_A).$$
\end{lem}

A (\textit{weakly}) \textit{monotone function} is a function $f:(X,\leqslant_X)\rightarrow(Y,\leqslant'_Y)$ between preordered sets satisfying $f(x)\leqslant'_Yf(y)$ whenever $x\leqslant_Xy$.

\subsection{Topology and order}\label{subsec:preordered.spaces}

We are interested in sets which are at once topologized and preordered, such as the circle $\mathbb{S}^1$ of states discussed in \S\ref{sec:intro}.
Recall that a \textit{preordered space} is a preordered set whose underlying set comes equipped with a topology - bearing no particular relationship with the order.

\begin{eg}\label{eg:customer.support}
	Let $n$ be a positive integer and $\pi_i$ be the $i$th projection map $\mathbb{I}^n\rightarrow\mathbb{I}$.
	Define $\leqslant_{\partial\mathbb{I}^n}$ to be the partial order on the subspace
	$$\partial\mathbb{I}^n=\{x\in\mathbb{I}^n\;|\;\prod_{i=1}^n\pi_i(x)(1-\pi_i(x))=0\}\subset\mathbb{I}^n$$
	with smallest graph satisfying $x\leqslant_{\partial\mathbb{I}^n}y$ whenever (i) $\pi_i(x)\leq\pi_i(y)$ for all integers $1\leq i\leq n$ and (ii) $\pi_j(x)=\pi_j(y)\in\{0,1\}$ for some integer $1\leq j\leq n$.

	Each point $x\in\partial\mathbb{I}^n$ can encode the collective progress of, say, $n$ customers on a telephone service with $n-1$ operators.  
	Each coordinate of $x$ represents the progress of an individual customer; at each point in time, at least one customer has made no progress $(0)$ or complete progress $(1)$.
	The partial order $\leqslant_{\partial\mathbb{I}^n}$ describes the causal relationship between the possible states of the customer-operator system.  

	The homotopy type of $\partial\mathbb{I}^n\simeq\mathbb{S}^{n-1}$ captures obstructions to simultaneous telephone support for all customers.
	As we glue such spaces together to form state spaces of more complicated systems, we require a subtler homotopy theory to remember the order of such obstructions.  
	See \cite{fgr:ditop} for a homotopy theory of partially ordered spaces.
\end{eg}

Following \cite{nachbin:order}, we call a preordered space \textit{locally convex} if it admits a basis of open convex subsets.

\begin{eg}\label{eg:pospaces}
	A \textit{pospace} is a partially ordered space $(X,\leqslant_X)$ such that $\mathrm{graph}(\leqslant_X)$ is closed in the standard product topology on $X\times X$.  
	Pospaces are automatically Hausdorff by \cite[Proposition 2]{nachbin:order}.  
	Compact pospaces are locally convex by \cite[Theorem 5]{nachbin:order}.
	The reader can check that the preordered space in Example \ref{eg:customer.support} is a pospace and hence locally convex.
\end{eg}

\begin{eg}[A non-example]
	The preordered circle of \S\ref{sec:intro} is not locally convex because its only convex subsets are $\varnothing$ and $\mathbb{S}^1$.
\end{eg}

\begin{eg}[A partially ordered non-example]
	Let $\leqslant_{\mathrm{lex}}$ denote the ``lexicographic order'' on $\mathbb{R}\times\mathbb{R}$ defined by the rule $(s_1,s_2)\leqslant_{\mathrm{lex}}(t_1,t_2)$ if $s_1<t_1$ or $s_1=t_1$ and $s_2\leq t_2$.  
	The totally ordered space $(\mathbb{R}\times\mathbb{R},\leqslant_{\mathrm{lex}})$ is not locally convex because images of its non-empty convex subsets under projection $\mathbb{R}\times\mathbb{R}\rightarrow\mathbb{R}$ onto the second factor are singletons and $\mathbb{R}$. 
\end{eg}

Local convexity sometimes implies ``antisymmetry.''  
Recall that a \textit{partially ordered space} is a preordered space whose preorder is a partial order.

\begin{lem}\label{lem:antisymmetric.convexity}
	A $T_0$ preordered space is a partially ordered space if each of its points admits a local base of convex neighborhoods.
\end{lem}
\begin{proof}
	Consider a preordered space $(X,\leqslant_X)$ and suppose $x\leqslant_Xy\leqslant_Xy$ for some $x,y\in X$.
	Then $x,y$ share the same convex neighborhoods and they hence share a common local base if they both admit a local base of convex neighborhoods.  
	Thus $x=y$ if $X$ is also $T_0$.
\end{proof}

Recall that a \textit{monotone map} $f:(X,\leqslant_X)\rightarrow(Y,\leqslant_Y)$ is a (weakly) monotone and continuous function between preordered spaces.

\begin{eg}
	Continuing Example \ref{eg:customer.support}, $x\leqslant_{\partial\mathbb{I}^n}y$ if and only if there is a monotone map $(\mathbb{I},\leqslant_{\mathbb{I}})\rightarrow(\partial\mathbb{I}^n,\leqslant_{\partial\mathbb{I}^n})$ sending $0$ to $x$ and $1$ to $y$.
	Following \cite{fgr:ditop}, we can think of such paths as all possible evolutions of our understaffed telephone support service.
\end{eg}

\section{Locally preordered spaces}\label{sec:streams}
We introduce the category of \textit{streams}, spaces equipped with some coherent preordering of their open subsets.
In \S\ref{subsec:definitions}, we introduce some basic definitions and examples.
We then construct limits and colimits in \S\ref{subsec:cosheafafications}.  

\subsection{Streams: basic definitions and examples}\label{subsec:definitions}

Given a state space $X$, we can write $x\leqslant_Uy$ whenever a machine can evolve from $x$ to $y$ when restricted to an open subset $U\subset X$ of states.  
Borrowing a metaphor from \cite{grandis:d}, we call such endowed spaces \textit{streams}.  

\begin{defn}\label{defn:streams}
	A \textit{circulation} $\leqslant$ \textit{on a space $X$} is a function assigning to each open subset $U\subset X$ a preorder $\leqslant_U$ on $U$ such that for each collection $\mathcal{O}$ of open subsets of $X$,
	\begin{equation}\label{eqn:cosheaf}
		\leqslant_{\bigcup\mathcal{O}}\;=\;\bigvee_{U\in\mathcal{O}}\leqslant_U.
	\end{equation}

	A \textit{stream} $(X,\leqslant)$ is a space $X$ equipped with a circulation $\leqslant$ on it.  
	We say that the preordered space $(X,\leqslant_X)$ \textit{underlies} a stream $(X,\leqslant)$.
\end{defn}

\begin{eg}
	For a preordered set $(X,\leqslant_X)$, the function sending $\varnothing$ to the unique preorder on $\varnothing$ and sending $X$ to $\leqslant_X$ defines a circulation on $(X,\{\varnothing,X\})$. 
\end{eg}

As another trivial example, we can form the ``weakest'' circulation on a given space.

\begin{defn}
	The \textit{trivial} circulation $\mathrm{id}$ on a space $X$ is the circulation assigning to each open subset $U\subset X$ the trivial preorder $\mathrm{id}_U$.
\end{defn}

We can also form the ``strongest'' circulation on a given space.

\begin{lem}\label{lem:pointwise}
	The pointwise application of $\bigvee$ to a non-empty set of circulations on a space is a circulation.
\end{lem}
\begin{proof}
	The operation $\bigvee$ is associative and commutative.
\end{proof}

The following examples suggest a relationship between connectivity and non-trivial circulation implicit in (\ref{eqn:cosheaf}).

\begin{eg}\label{eg:connect.circulations}
	Every space admits a circulation $\sim$ defined by the rule ``$x\!\sim_U\!y$ if $x,y$ occupy a common compact Hausdorff and connected subspace of $U$ whose closed subsets are locally connected.''  
\end{eg}

\begin{eg}\label{eg:specialization}
	Consider the \textit{specialization preorder} $\leqslant_X$ on a space $X$ defined by $x\leqslant_Xy$ if all neighborhoods of $x$ contain $y$.
	The rule 
	$$(\leqslant_X)_{\restriction(-)}:U\mapsto(\leqslant_X)_{\restriction U}$$
	defines a circulation because $x\leqslant_Uy$ whenever $x\leqslant_{\bigcup\mathcal{O}}y$, for each set $\mathcal{O}$ of open subsets of $X$ and each neighborhood $U\in\mathcal{O}$ of $x$ in $X$.
\end{eg}

\begin{eg}[Stream from d-space]\label{eg:d-spaces}
	Consider a \textit{d-space} $(X,dX)$ as defined in \cite{grandis:d}: a space $X$ and a set $dX$ on $X$ containing all constant paths and closed under both concatenation and precomposition with monotone maps $\mathbb{I}\rightarrow\mathbb{I}$. 
	The rule ``$x\leqslant^{dX}_Uy$ when $dX$ contains a path in $U$ from $x$ to $y$'' defines a circulation because every path in $U$ is the concatenation of paths whose images each lie inside elements of a given open cover of $U$.
\end{eg}

A restriction of machine states results in a restriction of machine behavior.

\begin{lem}\label{lem:precirculations}
	For all open subsets $U\subset V$ of a stream $(X,\leqslant)$, 
	$$\graph{\leqslant_U}\subset\graph{\leqslant_V}.$$
\end{lem}
\begin{proof}
	Observe
	$$\graph{\leqslant_U}\subset\graph{\bigvee\{\leqslant_U,\leqslant_V\}}=\graph{\leqslant_{U\cup V}}=\graph{\leqslant_V}.$$
\end{proof}

For a set $\mathcal{O}$ of open subsets of a stream $(X,\leqslant)$, each machine step 
$$x\leqslant_{\bigcup\mathcal{O}}y$$
breaks down into some sequence of smaller steps
$$x=x_0\leqslant_{U_1}x_1\leqslant_{U_2}\ldots\leqslant_{U_n}x_n=y$$
for some integer $n\geq 0$ and some $U_1,\ldots,U_n\in\mathcal{O}$ by (\ref{eqn:cosheaf}).  
When $\mathcal{O}$ contains exactly two sets, the following lemma asserts that we can take $U_1,\ldots,U_n$ to alternate and that consequently $x_1,\ldots,x_{n-1}\in\bigcap\mathcal{O}$.

\begin{lem}\label{lem:cosheaf}
	Consider a stream $(X,\leqslant)$.  
	For every pair of open subsets $U_1,U_{-1}\subset X$ and every pair $x,y\in X$ satisfying $x\leqslant_{U_1\cup U_{-1}}y$, there exist $n\in\{0,1,\ldots\}$, $m\in\{0,1\}$, and a sequence
	$$x=x_0\leqslant_{U_{(-1)^{m}}}x_1\leqslant_{U_{(-1)^{m+1}}}\cdots\leqslant_{U_{(-1)^{m+n}}}x_n=y.$$
\end{lem}
\begin{proof}
	For some minimal positive integer $n$, there exists a sequence of points $x=x_0,x_1,\cdots,x_n=y\in X$ such that for each $0<i\leq n$, either $x_{i-1}\leqslant_{U_1}x_i$ or $x_{i-1}\leqslant_{U_{-1}}x_i$.
	For all $0<j<n$, neither $x_{j-1}\leqslant_{U_1}x_j\leqslant_{U_1}x_{j+1}$ nor $x_{j-1}\leqslant_{U_{-1}}x_j\leqslant_{U_{-1}}x_{j+1}$ because $n$ is minimal and $\leqslant_{U_1},\leqslant_{U_{-1}}$ are transitive.  
\end{proof}

We can apply the lemma, for example, to extract information about underlying preordered spaces.

\begin{lem}\label{lem:connectedness}
	The bounded intervals of the underlying preordered space of every stream have connected closures.
\end{lem}
\begin{proof}
	Consider a stream $(X,\leqslant)$ and a non-empty bounded interval	
	$$I=\;\upper{X}{x}\;\cap\downer{X}{y}$$
	in $(X,\leqslant_X)$, for some $x,y\in X$.  
	Let $\bar{I}$ denote the closure in $X$ of $I$.  
	Consider some open subsets $U,V\subset X$ such that $\bar{I}\subset U\cup V$, $U\cap\bar{I}\neq\varnothing$, and $V\cap\bar{I}\neq\varnothing$.
	We can assume $x\in U$ without loss of generality and consider some $v\in V\cap I$.

	In the case $v\in U$, $v\in U\cap V\cap\bar{I}$.	
	Suppose $v\notin U$.  
	There is a sequence
	$$x=x_0\leqslant_Ux_1\leqslant_{V\cup(X\setminus\bar{I})}\cdots\leqslant_Ux_{n-1}\leqslant_{V\cup(X\setminus\bar{I})}x_n=v$$
	for some integer $n>0$ by Lemma \ref{lem:cosheaf}.  
	We conclude $x\leqslant_Xx_1\leqslant_Xv\leqslant_Xy$ from Lemma \ref{lem:precirculations} and $x_1\in U\cap (V\cup(X\setminus\bar{I}))$.  
	Hence $x_1\in U\cap (V\cup(X\setminus\bar{I}))\cap I\subset U\cap V\cap\bar{I}$.  
	In either case, $U\cap V\cap\bar{I}\neq\varnothing$.
\end{proof}

\begin{eg}\label{eg:trivial}
	A circulation on a totally disconnected space must be trivial by Lemma \ref{lem:connectedness}.  
	As likewise noted in \cite{bw:models}, ``atlases'' of partial orders are trivial on discrete spaces.
\end{eg}

Stream maps preserve all structure in sight.

\begin{defn}
	Given two streams $(X,\leqslant)$ and $(Y,\leqslant')$, a \textit{stream map $f:(X,\leqslant)\rightarrow(Y,\leqslant')$} is a continuous function $f:X\rightarrow Y$ satisfying $f(x)\leqslant_Uf(y)$ whenever $x\leqslant_{f^{-1}U}y$, for all open subsets $U\subset Y$. 
\end{defn}

Let $\STREAMS$ denote the category of streams and stream maps.

\begin{eg}[d-space from stream]\label{eg:dipaths}
	Continuing Example \ref{eg:d-spaces}, we can take $d\mathbb{I}$ to be the set of monotone paths $\mathbb{I}\rightarrow\mathbb{I}$ and we can associate to each stream $(Y,\leqslant')$ the d-space $(Y,\STREAMS((\mathbb{I},\leqslant^{d\mathbb{I}}),(Y,\leqslant')))$.
	Our constructions from streams to d-spaces and vice versa extend to an adjunction between $\STREAMS$ and a category $d\SPACES$ of d-spaces defined in \cite{grandis:d}. 
\end{eg}

\subsection{Cosheafifications of precirculations}\label{subsec:cosheafafications}

Presheaves are to sheaves what \textit{precirculations} are to circulations: ``functors'' which need no longer satisfy the ``sheaf condition'' (\ref{eqn:cosheaf}).
We can construct with ease a precirculation satisfying some prescribed universal property, and then apply \textit{cosheafification} to obtain some desired universal stream. 
\begin{defn}
	A \textit{precirculation $\leqslant$ on a space $X$} is a function assigning to each open subset $U\subset X$ a preorder $\leqslant_U$ on $U$ such that $\graph{\leqslant_U}\subset\graph{\leqslant_V}$ whenever $U\subset V$.
\end{defn}

Precirculations generalize circulations by Lemma \ref{lem:precirculations} and satisfy half of the ``cosheaf'' condition (\ref{eqn:cosheaf}).

\begin{lem}\label{lem:half.cosheaf}
	For every precirculation $\leqslant$ on a space $X$ and every set $\mathcal{O}$ of open subsets of $X$,
	$$\graph{\!\bigvee_{U\in\mathcal{O}}\!\!\leqslant_U}\subset\graph{\leqslant_{\bigcup\mathcal{O}}}.$$
\end{lem}
\begin{proof}
	Immediate from Lemma \ref{lem:transitive.closure}.
\end{proof}

Our generalizations of circulations straightforwardly ``pullback'' and ``pushforward'' along continuous functions.

\begin{defn}
	Consider a map $f:X\rightarrow Y$ of spaces and let $\mathcal{T}_Y$ denote the topology of $Y$.  
	For each precirculation $\leqslant$ on $Y$, the \textit{pullback $\leqslant^{f^*}$ of $\leqslant$ along $f$} is the precirculation on $X$ defined by
	$$\graph{\leqslant^{f^*}_U}=\!\!\!\!\!\bigcap_{f(U)\subset V\in\mathcal{T}_Y}\!\!\!\!\!(f\times f)^{-1}(\graph{\leqslant_V})\cap(U\times U).$$

	For each precirculation $\leqslant'$ on $X$, the \textit{pushforward $(\leqslant')^{f_*}$ of $\leqslant'$ along $f$} is the precirculation on $Y$ assigning to each open subset $U\subset Y$ the transitive-reflexive closure of the relation on $U$ with graph 
	$$(f\times f)(\graph{\leqslant'_{f^{-1}U}}).$$
\end{defn}

\begin{eg}
	Consider a continuous function $f:X\rightarrow Y$.  
	Then $\mathrm{id}^{f_*}=\mathrm{id}$ but $\mathrm{id}^{f^*}=\mathrm{id}$ if and only if $f$ is injective.
\end{eg}

Pushfowards and certain pullbacks preserve circulations.  

\begin{lem}\label{lem:pushforwards}
	For each continuous function $f:X\rightarrow Y$ of spaces and each circulation $\leqslant$ on $X$, $\leqslant^{f_*}$ is a circulation.
\end{lem}
\begin{proof}
	Consider a family $\mathcal{O}$ of open subsets of $Y$. 
	Then
	\begin{eqnarray*}
		(f\times f)(\graph{\leqslant_{f^{-1}\bigcup\mathcal{O}}}) 
		&=& (f\times f)(\graph{\bigvee_{U\in\mathcal{O}}\leqslant_{f^{-1}U}}),\\
		(f\times f)(\bigcup_{U\in\mathcal{O}}\graph{\leqslant_{f^{-1}U}})
		&=& \bigcup_{U\in\mathcal{O}}(f\times f)(\graph{\leqslant_{f^{-1}U}})\\
		&\subset& \graph{\bigvee_{U\in\mathcal{O}}\leqslant^{f_*}_U},
	\end{eqnarray*}
	and therefore $\graph{\leqslant^{f_*}_{\bigcup\mathcal{O}}}\subset\graph{\bigvee_{U\in\mathcal{O}}\leqslant^{f_*}_U}$ by Lemma \ref{lem:transitive.closure}.  
	The result follows from Lemma \ref{lem:half.cosheaf}.
\end{proof}

\begin{eg}\label{eg:quotients}
	Continuing Example \ref{eg:d-spaces}, we can model the states $\mathbb{S}^1$ of \S\ref{sec:intro} as the stream $(\mathbb{S}^1,(\leqslant^{d\mathbb{I}})^{q_*})$, where $q:\mathbb{I}\rightarrow\mathbb{S}$ denotes the quotient map identifying endpoints.  
\end{eg}

\begin{eg}\label{eg:pullbacks}
	Pullbacks along open maps similarly preserve circulations.
\end{eg}

\begin{eg}[Pullback not preserving circulation]\label{eg:pullback.pathology}
	Let $\Delta$ be the diagonal map $\mathbb{I}\rightarrow\mathbb{I}^2$ and let $d\mathbb{I}^2$ denote the set of all paths on $\mathbb{I}^2$ which are (weakly) monotone in each coordinate and whose images intersect $\Delta(0,1)$ only when those images are singletons.  
	Continuing Example \ref{eg:d-spaces}, let $\leqslant=\leqslant^{d\mathbb{I}^2}$.
	Although $0\leqslant^{\Delta^*}_{\mathbb{I}}\!\!1$ because every neighborhood in $\mathbb{I}^2$ of $\Delta(\mathbb{I})$ contains a path in $d\mathbb{I}^2$ from $\Delta(0)$ to $\Delta(1)$, 
	$$\leqslant^{\Delta^*}_{\mathbb{I}}\!\![0]\cap(\leqslant^{\Delta^*}_{\mathbb{I}})^{-1}[1]=\{0,1\}$$
	has disconnected closure and hence $\leqslant^{\Delta^*}_{\mathbb{I}}$ cannot be a circulation by Lemma \ref{lem:connectedness}.
\end{eg}

Pushforwards satisfy a universal property best articulated in the language of categorical topology.
Consider a functor $F:\GENERIC_1\rightarrow\GENERIC_2$ - our motivating example is the forgetful functor $T:\STREAMS\rightarrow\SPACES$.  
Consider also a functor $D:\DIAGRAM\rightarrow\GENERIC_1$ and some cocones on $D$ and $FD$ as follows:
$$\lambda=(\lambda_x:Dx\rightarrow c)_{x\in\mathrm{ob}\;\DIAGRAM},\quad\bar\lambda=(\bar{\lambda}_x:FDx\rightarrow\bar{c})_{x\in\mathrm{ob}\;\DIAGRAM}.$$
Dualizing the definition of \textit{initial structures} in \cite{borceux:2}, we call $\lambda$ a \textit{final structure for $(F,D,\bar\lambda)$ to $c$} if $F\lambda=\bar\lambda$ and for each cocone $\lambda'$ from $D$ to a $\GENERIC_1$-object $c'$, every $\GENERIC_2$-map $Fc\rightarrow Fc'$ by which $F\lambda'$ factors through $F\lambda$ has a unique preimage under $F$ by which $\lambda'$ factors through $\lambda$.

The following lemma asserts that final structures always exist for the case $F=T$.  
In stating and proving the lemma, we adopt the following conventions.  
We reuse the symbol $\leqslant$ to signify different circulations.  
Also, we write $\bigvee$ for the pointwise application of $\bigvee$ to precirculations on a fixed space, defining $\bigvee\varnothing=\mathrm{id}$.

\begin{lem}\label{lem:final}
	For each diagram $D:\DIAGRAM\rightarrow\STREAMS$ and each cocone
	$$\lambda=(\lambda_x:TDx\rightarrow X)_{x\in\mathrm{ob}\;\DIAGRAM}$$
	on $TD$, $(X,\bigvee_{x\in\mathrm{ob}\;\DIAGRAM}\leqslant^{(\lambda_x)_*})$ is a stream to which there is a final structure for $(T,D,\lambda)$.
\end{lem}
\begin{proof}
	The precirculation $\leqslant'=\bigvee_{x\in\mathrm{ob}\;\DIAGRAM}\leqslant^{(\lambda_x)_*}$ is a circulation by Lemmas \ref{lem:pointwise} and \ref{lem:pushforwards}.  
	For all $x\in\mathrm{ob}\;\DIAGRAM$, the function $\lambda_x$ defines a stream map $Dx\rightarrow(X,\leqslant')$ because for all open subsets $U\subset X$,
	$$(\lambda_x\times\lambda_x)(\graph{\leqslant_{\lambda_x^{-1}U}})\subset\graph{\leqslant'_U}.$$

	Consider a cocone $(\lambda'_x:Dx\rightarrow(Y,\leqslant''))_{x\in\mathrm{ob}\;\DIAGRAM}$ on $D$ and a continuous function $\eta:X\rightarrow Y$ such that $\eta\lambda_x=\lambda'_x$ for each $x\in\mathrm{ob}\;\DIAGRAM$.
	For all open subsets $V\subset Y$ and all $x\in\mathrm{ob}\;\DIAGRAM$, 
	\begin{eqnarray*}
		(\eta\times\eta)(\lambda_x\times\lambda_x)(\graph{\leqslant_{\lambda_x^{-1}(\eta^{-1}V)}})
		&=&(\lambda'_x\times\lambda'_x)(\graph{\leqslant_{(\lambda'_x)^{-1}V}})\\
		&\subset&\graph{\leqslant''_V},
	\end{eqnarray*}
	hence $(\eta\times\eta)(\graph{\leqslant^{(\lambda_x)_*}_{\eta^{-1}V}}\subset\graph{\leqslant''_V})$ by Lemma \ref{lem:transitive.closure}, and hence
	$$(\eta\times\eta)(\graph{\leqslant'_{\eta^{-1}V}})\subset\graph{\leqslant''_V}$$
	by Lemma \ref{lem:transitive.closure} again.  
	Equivalently, $\eta$ defines a stream map $(X,\leqslant')\rightarrow(Y,\leqslant'')$.  
\end{proof}

In the language of \cite{borceux:2}, a functor $F:\GENERIC_1\rightarrow\GENERIC_2$ is \textit{topological} if there always exists an initial structure for the appropriate set $(F,D,\lambda)$ of data.  
As a consequence of Lemma \ref{lem:final} and \cite[Proposition 7.3.6, Proposition 7.3.11]{borceux:2}, the functor $T$ is topological and thus creates limits and colimits by \cite[Proposition 7.3.8]{borceux:2}.  
We form initial structures, such as limits, by approximating precirculations with circulations in a manner analogous to sheafification.

\begin{defn}
	The \textit{cosheafification $\leqslant^!$} of a precirculation $\leqslant$ on a space $X$ is the pointwise application of $\bigvee$ to all circulations $\leqslant'$ on $X$ satisfying $\graph{\leqslant'_U}\subset\graph{\leqslant_U}$ for each open subset $U\subset X$.
\end{defn}

\begin{eg}[Stream from locally partially ordered space]\label{eg:atlases}
	A \textit{local partial order} $\leqslant$ on a space $X$, defined in \cite{fgr:ditop}, amounts to a precirculation on $X$ such that for some open cover $\mathcal{O}$ of $X$, $\leqslant$ assumes partial orders on $\mathcal{O}$ and for all open subsets $U\subset X$ and all $V\in\mathcal{O}$, $\leqslant_{U\cap V}=(\leqslant_V)_{\restriction U\cap V}$.
	Two local partial orders are \textit{equivalent} if they coincide on an open basis.
	
	The cosheafification of a local partial order $\leqslant$ only depends upon the equivalence class $[\leqslant]$ of $\leqslant$.
	We thereby obtain a functor from a category in \cite{bw:models} of \textit{locally partially ordered spaces} $(X,[\leqslant])$ to $\STREAMS$.
	Images of locally partially ordered spaces under this functor include the ``quotient stream'' constructed in Example \ref{eg:quotients}.
\end{eg}

\begin{eg}[Product streams]\label{eg:products}
	Consider a family $\{(X_i,\leqslant^i)\}_{i\in\mathcal{I}}$ of streams and let $X=\prod_iX_i$.  
	For each $j\in\mathcal{I}$, the projection $\pi_j:X\rightarrow X_j$ onto the $j$th factor is an open map.  
	We can therefore define a precirculation $\leqslant$ on $X$ by the rule
	$$\graph{\leqslant_U}=\graph{\prod_{i\in\mathcal{I}}\leqslant^i_{\pi_i(U)}}\cap(U\times U).$$
	This precirculation is almost never a circulation on $X$.  
	We leave it to the reader to check that $(X,\leqslant^!)$ is the $\STREAMS$-product $\prod_{i\in\mathcal{I}}(X_i,\leqslant^i)$.
\end{eg}

The general failure of pullbacks along inclusions to preserve circulations, as demonstrated in Example \ref{eg:pullback.pathology}, reflects some of the pathology that \textit{substreams} can exhibit.

\begin{defn}
	For a stream $(X,\leqslant)$ and a subspace $A\subset X$,
	$$(A,(\leqslant^{(A\hookrightarrow X)^*})^!)$$
	is a \textit{substream} of $(X,\leqslant)$ and we write $\leqslant_{\restriction A}$ for $(\leqslant^{(A\hookrightarrow X)^*})^!_A$.
	A \textit{stream inclusion} is an inclusion function to a stream $(X,\leqslant)$ from a substream of $(X,\leqslant)$.
\end{defn}

For example, an open substream of a stream $(X,\leqslant)$ is just an open subspace of $X$ equipped with a suitable restriction of $\leqslant$.  
We give a concrete description of certain closed substreams below.

\begin{eg}\label{eg:convex}
	Let $(X,\leqslant)$ be a stream and let $A$ be a closed and convex subspace of $(X,\leqslant_X)$.
	The precirculation $\leqslant'$ on $A$ sending $A\cap U$ to $(\leqslant_U)_{\restriction A\cap U}$ for each open subset $U\subset X$ is well-defined - for open subsets $V,W\subset X$ such that $A\cap V=A\cap W$, 
	\begin{eqnarray*}
	(\leqslant_V)_{\restriction A\cap V}
	&=& (\leqslant_V)_{\restriction A\cap V\cap W}\\
	&=&(\leqslant_{V\cap W})_{\restriction A\cap V\cap W}\\
	&=&(\leqslant_W)_{\restriction A\cap V\cap W}\\
	&=&(\leqslant_W)_{\restriction A\cap W}
	\end{eqnarray*}
	by an application of Lemma \ref{lem:convex} given later in \S\ref{sec:nature}.
	It is straightforward to check that $\leqslant'$ is a circulation because $\leqslant$ is a circulation and $A$ is convex in $(X,\leqslant_X)$.  
	It follows that $(A,\leqslant')$ is a substream of $(X,\leqslant)$.	
\end{eg}

\begin{eg}[One-point compactifications]\label{eg:one.point}
	Let $\hat{X}=X\cup\{\infty\}$ be the one-point compactification of a non-compact space $X$.
	For each circulation $\leqslant$ on $X$, define a precirculation $\hat{\leqslant}$ on $\hat{X}$ as sending each open subset $U\subset\hat{X}$ to the preorder on $U$ with smallest graph containing
	$$(U\times U)\cap((\{\infty\}\times X)\cup(X\times\{\infty\})\cup\graph{\leqslant_{X\cap U}}).$$
	
	The reader can check that for all locally compact Hausdorff, non-compact streams $(X,\leqslant)$, $(\hat{X},\hat{\leqslant}\;\!^!)$ is ``final'' among all compact Hausdorff streams in which $(X,\leqslant)$ lies as a dense substream.  
\end{eg}

As a formal consequence of our definitions, stream inclusions define ``initial structures.''
We include proofs of the next two lemmas for completeness.

\begin{lem}\label{lem:inclusions}
	Stream inclusions are stream maps.  
	For each diagram
	$$\xymatrix{
	(X,\leqslant)\ar@{.>}[dr]\ar[rr]^f & & (Y,\leqslant')\\
	& (A,\leqslant^A)\ar[ur]_\iota
	}$$
	of solid arrows denoting stream maps, where $\iota$ denotes stream inclusion and $f(X)\subset A$, there exists a unique dotted stream map making the diagram commute. 
\end{lem}
\begin{proof}
	Let $\DIAGRAM$ be the discrete category of streams $(A,\leqslant'')$ from which inclusion defines a stream map to $(Y,\leqslant')$.
	Consider a circulation $\leqslant'''$ on $A$.  
	Then $(A,\leqslant''')\in\mathrm{ob}\;\DIAGRAM$ if and only if, for each open subset $V\subset Y$,
	$$\graph{\leqslant'''_U}\subset(U\times U)\cap\graph{\leqslant'_V}$$
	for $U=A\cap V$ and hence for all open subsets $U\subset A$ sitting inside $V$ by Lemma \ref{lem:precirculations}.
	Thus $(A,\leqslant''')\in\mathrm{ob}\;\DIAGRAM$ if and only if, for each open subset $U\subset A$,
	\begin{equation}\label{eqn:inclusion.inclusion}
	\graph{\leqslant'''_U}\subset\graph{(\leqslant')^{(A\hookrightarrow Y)^*}_U}.
	\end{equation}

	Let $D$ be the inclusion functor $\DIAGRAM\hookrightarrow\STREAMS$, and let $\bar\lambda$ be the cocone on $TD$ to $A$ defined by identity maps.  
	The final structure for $(T,D,\bar\lambda)$ is a cocone $\lambda$ to $(A,\leqslant^A)$ by Lemma \ref{lem:final} because
	\begin{eqnarray*}
		\leqslant^A
		&=& \bigvee_{(A,\leqslant'')\in\mathrm{ob}\;\DIAGRAM}\leqslant''\\
		&=& \bigvee_{(A,\leqslant'')\in\mathrm{ob}\;\DIAGRAM}(\leqslant'')^{(\mathrm{id}_A)_*},
	\end{eqnarray*}
	the first line following from (\ref{eqn:inclusion.inclusion}).
	One consequence is that inclusion defines a stream map $(A,\leqslant^A)\hookrightarrow(Y,\leqslant')$.  

	The continuous function $T(f)$ factors through its corestriction $\bar{g}:X\rightarrow A$ and $T(\iota)$.  
	Let $E:\{*\}\rightarrow\STREAMS$ be the discrete diagram sending $*$ to $(X,\leqslant)$, and let $\bar\mu$ be the cocone on $TE$ to $A$ defined by $\bar{g}$.
	There exists a final structure $\mu$ for $(T,E,\bar\mu)$ to $(A,\leqslant^{\bar{g}_*})$ by Lemma \ref{lem:final}.
	In particular, $\bar{g}$ defines a stream map $g:(X,\leqslant)\rightarrow(A,\leqslant^{\bar{g}_*})$.  
	Let $\mu'$ denote the cocone on $E$ to $(Y,\leqslant')$ defined by $f$.
	The stream $(A,\leqslant^{\bar{g}_*})\in\DIAGRAM$ because $\mu'$ factors through $\mu$.
	The stream map 
	$$\lambda_{(A,\leqslant^{\bar{g}_*})}:(A,\leqslant^{\bar{g}_*})\rightarrow(A,\leqslant^A)$$
	composes with $g$ to define our dotted arrow.  
	Uniqueness follows from the uniqueness of the dotted arrow making the diagram of underlying sets commute.
\end{proof}

A consequence is that a substream of a substream of a stream $(X,\leqslant)$ is a substream of $(X,\leqslant)$.

\begin{lem}\label{lem:inclusion.composites}
	The composite of stream inclusions is a stream inclusion.
\end{lem}
\begin{proof}
	Consider three stream inclusions	
	$$i:(B,\leqslant'')\hookrightarrow(C,\leqslant'),\;j:(C,\leqslant')\hookrightarrow(D,\leqslant),\;k:(B,\leqslant''')\hookrightarrow(D,\leqslant).$$
	Then $\id_B$ defines a stream map $(B,\leqslant'')\rightarrow(B,\leqslant''')$ by Lemma \ref{lem:inclusions} applied to the case $A=B$ and $f=ji$.  
	Inclusion defines a stream map $i':(B,\leqslant''')\rightarrow(C,\leqslant')$ by Lemma \ref{lem:inclusions} applied to the case $A=C$ and $f=k$.
	Thus $\id_B$ defines a stream map $(B,\leqslant''')\rightarrow(B,\leqslant'')$ by Lemma \ref{lem:inclusions} applied to the case $A=B$ and $f=i'$.
\end{proof}

Circulations of substreams sometimes obey a generalization of (\ref{eqn:cosheaf}).

\begin{lem}\label{lem:pseudo.circulations}
	Consider a stream $(X,\leqslant)$.
	For each family $\mathcal{N}$ of subsets of $X$ containing a neighborhood in $X$ for each point in its union,
	\begin{equation}\label{eqn:pseudo.circulations}
		\leqslant_{\bigcup\mathcal{N}}\;=\;\bigvee_{A\in\mathcal{N}}\leqslant_{\restriction A}\;=\;\bigvee_{A\in\mathcal{N}}(\leqslant_{\bigcup\mathcal{N}})_{\restriction A}.
	\end{equation}
\end{lem}
\begin{proof}
	Let $A^\circ$ denote the interior in $X$ of a subset $A\subset X$.  
	Then
	\begin{eqnarray*}
	\graph{\leqslant_{A^\circ}}
	&=& \graph{\leqslant_{\restriction A^\circ}}\\
	&\subset& \graph{\leqslant_{\restriction A}}
	\end{eqnarray*}
	for all $A\in\mathcal{N}$, the first line due to the fact that circulations of open substreams of $(X,\leqslant)$ are just restrictions of $\leqslant$ and the second line due to an application of Lemma \ref{lem:inclusions}.
	Thus in (\ref{eqn:pseudo.circulations}), the graph of the leftmost preorder lies inside the graph of the middle preorder because $\leqslant_{\bigcup\mathcal{N}}=\bigvee_{A\in\mathcal{N}}\leqslant_{A^\circ}$.
	For each $A\in\mathcal{N}$,
	\begin{eqnarray*}
	\graph{\leqslant_{\restriction A}}
	&=& \graph{(\leqslant^{(A\hookrightarrow X)^*})^!_{A}}\\
	&\subset& \graph{(\leqslant^{(A\hookrightarrow X)^*})_{A}}\\
	&\subset& \graph{(\leqslant_{\bigcup\mathcal{N}})_{\restriction A}}
	\end{eqnarray*}
	and therefore in (\ref{eqn:pseudo.circulations}), the graph of the middle preorder lies inside the graph of the rightmost preorder.
	Finally, the graph of the rightmost preorder lies inside the graph of the leftmost preorder in (\ref{eqn:pseudo.circulations}), completing the chain of inclusions demonstrating (\ref{eqn:pseudo.circulations}).
\end{proof}

\section{Streams in nature}\label{sec:nature}

Although a general stream appears to encode a frightening amount of order-theoretic information, we can identify streams in nature which are determined by their underlying preordered spaces.
We can start by identifying when ``global sections'' determine ``local sections.''

\begin{lem}\label{lem:convex}
	For a stream $(X,\leqslant)$ and a convex subset $A$ of $(X,\leqslant_X)$,
	$$(\leqslant_U)_{\restriction A}=(\leqslant_X)_{\restriction A}$$
	for each open neighborhood $U$ of the closure of $A$ in $X$.
\end{lem}
\begin{proof}
	Consider some $x,y\in A$ satisfying $x\leqslant_Xy$ and consider an open neighborhood $U$ of the closure $\bar{A}$ of $A$ in $X$.  
	There exists a sequence
	$$x=x_0\leqslant_Ux_1\leqslant_{X\setminus\bar{A}}\ldots\leqslant_Ux_n=y$$
	by Lemma \ref{lem:cosheaf}.  
	Then $n<2$ - otherwise $x_1\notin A$ because $x_1\in U\cap (X\setminus\bar{A})$ and $x\leqslant_X x_1\leqslant_Xy$ by Lemma \ref{lem:precirculations}, contradicting $A$ convex.  
	The result follows.
\end{proof}

Recall that a space $X$ is \textit{regular} if every closed subset $C\subset X$ and every point in $X\setminus C$ can be separated by disjoint neighborhoods. 

\begin{lem}\label{lem:uniquely.underlying}
	A regular Hausdorff preordered space underlies at most one stream if each of its points admits a local base of convex neighborhoods.`
\end{lem}
\begin{proof}
	Consider a regular Hausdorff stream $(X,\leqslant)$ whose points admit local bases of convex neighborhoods in $(X,\leqslant_X)$.
	Let $\mathcal{N}$ be the family of convex subsets of $(X,\leqslant_X)$ whose closures in $X$ lie inside an open subset $U\subset X$.  
	Then
	$$\leqslant_U\;=\bigvee_{A\in\mathcal{N}}(\leqslant_U)_{\restriction A}=\bigvee_{A\in\mathcal{N}}(\leqslant_X)_{\restriction A}.$$ 
	The first equality follows from Lemma \ref{lem:pseudo.circulations} because each point in a regular Hausdorff space has a local base of closed neighborhoods.
	The second equality follows from Lemma \ref{lem:convex}.
\end{proof}

\begin{eg}[Distinct d-spaces forming same stream]\label{eg:d.beaucoup}
	Let $s\mathbb{R}^2$ be the set of all concatenations of all paths on $\mathbb{R}^2$ which are constant in one coordinate and (weakly) monotone in the other coordinate.  
	Let $d\mathbb{R}^2$ be the set of all paths on $\mathbb{I}^2$ which are (weakly) monotone in both coordinates.
	Continuing Example \ref{eg:d-spaces}, $\leqslant^{s\mathbb{R}^2}=\leqslant^{d\mathbb{R}^2}$ by Lemma \ref{lem:uniquely.underlying}, even though $s\mathbb{R}^2\neq d\mathbb{R}^2$.
\end{eg}

We now give a criterion for a preordered space to underlie a stream.  

\begin{lem}\label{lem:underlies}
	A locally convex, compact Hausdorff preordered space underlies a stream if its bounded intervals are closed and connected.
\end{lem}
\begin{proof}
	Consider a locally convex, compact Hausdorff preordered space $(X,\leqslant_X)$ whose bounded intervals $[x,y]=\upper{X}{x}\;\cap\downer{X}{y}$ are closed and connected for all $x,y\in X$.  
	Define a precirculation $\leqslant'$ on $X$ by
	$$\leqslant'_U\;=\bigvee_{\varnothing\neq[x,y]\subset U}(\leqslant_X)_{\restriction\{x,y\}}$$

	Let $\BASIS$ denote the family of all open and convex subsets of $X$. 
	Consider a family $\mathcal{O}$ of open subsets of $X$ and a bounded interval $[a,b]\neq\varnothing$ contained in the union of $\mathcal{O}$.  
	There exists a finite subset $\mathcal{O}_F\subset\BASIS$ whose union also contains $[a,b]$ and whose elements each lie inside an element of $\mathcal{O}$ because $[a,b]$ is compact and $\BASIS$ is a basis.
	In order to show $a\;(\bigvee_{U\in\mathcal{O}}\!\leqslant'_U)\;b$ and thereby conclude $\leqslant'$ is a circulation from Lemma \ref{lem:half.cosheaf}, it suffices to show
	\begin{equation}\label{eqn:underlies}
		\graph{\leqslant'_{\bigcup\mathcal{O}_F}}\subset\graph{\bigvee_{U\in\mathcal{O}_F}\leqslant'_U}
	\end{equation}
	for such finite $\mathcal{O}_F\subset\BASIS$ by induction on the cardinality of $\mathcal{O}_F$.
	For then,
	$$(a,b)\in\graph{\leqslant'_{\bigcup\mathcal{O}_F}}\subset\graph{\!\!\bigvee_{U\in\mathcal{O}_F}\!\!\leqslant'_{U}}\subset\graph{\!\bigvee_{U\in\mathcal{O}}\!\leqslant'_U},$$
	the last containment following from the fact that $\leqslant'$ is a precirculation.
	
	As the base case $\mathcal{O}_F=\varnothing$ is vacuous, assume (\ref{eqn:underlies}) for all subsets $\mathcal{O}_F\subset\BASIS$ of cardinality $n-1$ and consider some bounded interval $[c,d]\neq\varnothing$ and some $U_1,\ldots,U_n\in\BASIS$ whose union contains $[c,d]$.  
	Heading towards a contradiction, assume $(c,d)\notin\graph{\bigvee_{i=1}^n\leqslant'_{U_i}}$.  
	We can assume $c\in U_1$ and $U_1\cap U_2\cap [c,d]$ contains some point $c'$ by reordering and $[c,d]$ connected.  
	
	The set $[c',d]\nsubseteq U_2\cup\ldots\cup U_{n}$ - otherwise $c\leqslant'_{U_1}c'\;(\bigvee_{i=2}^{n}\leqslant'_{U_i})\;d$ by our inductive hypothesis.
	Thus $[c',d]\cap(X\setminus (U_2\cup\ldots\cup U_n))$ contains some $c''$.  
	The set $[c'',d]\cap U_2$ contains some $c'''$ - otherwise $c\leqslant'_{U_1}c''\;(\leqslant'_{U_1}\!\!\!\!\vee\bigvee_{i=3}^n\leqslant'_{U_i})\;d$ by our inductive hypothesis.  
	However, $c',c'''\in U_2$ and $c''\notin U_2$, contradicting $U_2$ convex.
\end{proof}

\begin{eg}\label{eg:postreams}
	Continuing Example \ref{eg:pospaces}, a compact pospace underlies a stream if and only if its bounded intervals are connected by \cite[Propositions 1 and 2, Theorem 5]{nachbin:order} (or the more accessible \cite[Propositions VI-1.4 and VI-1.6, Corollary VI-1.9]{scott:lattices}) and Lemmas \ref{lem:connectedness} and \ref{lem:underlies}.  
	Examples of such pospaces include all connected, compact Hausdorff topological lattices by \cite[Proposition VI-5.15]{scott:lattices} and the pospace in Example \ref{eg:customer.support}. 
\end{eg}

\begin{eg}[A generalization]
	The hypothesis of local convexity in Lemma \ref{lem:underlies} is unnecessary.
	Consider a compact Hausdorff partially ordered space $(X,\leqslant_X)$ whose bounded intervals are closed and connected.  
	Let $\mathcal{M}$ denote the family of maximal chains in $(X,\leqslant_X)$ endowed with their subspace topologies.

	Consider an $A\in\mathcal{M}$.
	The space $A$ is closed in $X$ and hence compact Hausdorff by \cite[Proposition VI-5.1]{scott:lattices}.
	The pospace $(A,(\leqslant_X)_{\restriction A})$ is \textit{order-dense}, or equivalently its bounded intervals never contain exactly two points, by the maximality of $A$.
	Then $(A,(\leqslant_X)_{\restriction A})$ is a compact pospace whose intervals are connected and whose partial order is a total order, by \cite[Proposition VI-5.6]{scott:lattices}.  
	Therefore $(A,(\leqslant_X)_{\restriction A})$ underlies a stream $(A,\leqslant^A)$ by Lemma \ref{lem:underlies} and
	\begin{eqnarray*}
	\leqslant_X
	&=& \bigvee_{A\in\mathcal{M}}(\leqslant_X)_{\restriction A}\\
	&=& \bigvee_{A\in\mathcal{M}}(\leqslant_X)_{\restriction A}\vee\id_X\\
	&=& \bigvee_{A\in\mathcal{M}}(\leqslant^A)^{(A\hookrightarrow X)_*}_X\\
	&=& (\bigvee_{A\in\mathcal{M}}(\leqslant^A)^{(A\hookrightarrow X)_*})_X.
	\end{eqnarray*}
	We conclude that $(X,\leqslant_X)$ underlies $(X,\bigvee_{A\in\mathcal{M}}(\leqslant^A)^{(A\hookrightarrow X)_*})$.
\end{eg}

Let $\COMPACT$ denote the category of locally convex, compact Hausdorff partially ordered spaces whose bounded intervals are closed and connected, and monotone maps between them.

\begin{thm}\label{thm:embedding}
	There exists a full and concrete embedding
	\begin{equation}\label{eqn:embedding}
		\COMPACT\hookrightarrow\STREAMS
	\end{equation}
	sending each $\COMPACT$-object to a unique stream it underlies.  
	The image of (\ref{eqn:embedding}) contains all compact Hausdorff streams having locally convex underlying preordered spaces whose bounded intervals are closed.
\end{thm}
\begin{proof}
	Lemmas \ref{lem:uniquely.underlying} and \ref{lem:underlies} define the injective object class function $E$ of our desired embedding (\ref{eqn:embedding}).
	Consider some images $(X,\leqslant)$ and $(Y,\leqslant')$ of $E$ and a monotone map $f:(X,\leqslant_X)\rightarrow(Y,\leqslant'_Y)$.  
	Then
	$$(f\times f)(\graph{\leqslant_X})\subset\graph{\leqslant'_Y}$$
	or equivalently, $\leqslant'_Y=\bigvee\{\leqslant^{f_*}_Y,\leqslant'_Y\}$.  
	Therefore $\leqslant'=\bigvee\{\leqslant^{f_*},\leqslant'\}$ by Lemmas \ref{lem:pointwise}, \ref{lem:pushforwards}, and \ref{lem:uniquely.underlying}, and so $f$ defines a stream map $(X,\leqslant)\rightarrow(Y,\leqslant')$ by Lemma \ref{lem:final}.
	We conclude $E$ extends to a concrete functor, which is full because every stream map is a monotone map of underlying preordered spaces.
	The last statement of the theorem follows from Lemmas \ref{lem:antisymmetric.convexity} and \ref{lem:connectedness}.
\end{proof}

\begin{eg}\label{eg:local.embed}
	The functor defined in Example \ref{eg:atlases} restricts to an embedding from the category of those locally partially ordered spaces $(X,[\leqslant])$ admitting bases of open subsets $U$ whose points admit neighborhoods $(N,(\leqslant_U)_{\restriction N})$ in $\COMPACT$.
\end{eg}

\section{Compactly flowing streams}\label{sec:compactly.flowing}

We modify the category of streams in order to obtain a category of \textit{compactly flowing} streams whose limits and colimits are created by forgetting to the standard category of compactly generated spaces.  
We observe in this section that the category is convenient in the sense of \cite{steenrod:convenience} - it is Cartesian closed, it includes all streams of possible interest, and it is often closed under intuitive constructions.

Recall that a space $X$ is \textit{compactly generated} if it is weak Hausdorff (continuous images of compact Hausdorff spaces in $X$ are closed) and a subset $U\subset X$ is open in $X$ if $U\cap K$ is open in $K$ for all compact Hausdorff subspaces $K\subset X$. 
We can define when a circulation is also ``compactly generated.''  

\begin{defn}\label{defn:k}
	Fix a weak Hausdorff space $X$.
	Let $\mathcal{K}(X)$ denote the family of compact Hausdorff subsets of $X$.
	A circulation $\leqslant$ on $X$ is a \textit{k-circulation} if, for all open subspaces $U\subset X$,
	\begin{equation}\label{eqn:k.cosheaf}
		\leqslant_U\;=\bigvee_{L\in\mathcal{K}(U)}\leqslant_{\restriction L}.
	\end{equation}
	A stream $(X,\leqslant)$ is \textit{compactly flowing} if $X$ is compactly generated and $\leqslant$ is a k-circulation.
\end{defn}

Equivalently, a circulation $\leqslant$ on a weak Hausdorff space $X$ is a k-circulation if and only if, for each open subspace $U\subset X$,
\begin{equation}\label{eqn:alternate.k.circulation}
	\graph{\leqslant_U}=\bigcup_{L\in\mathcal{K}(U)}\graph{\leqslant_{\restriction L}}
\end{equation}
because $\mathcal{K}(U)$ is closed under finite unions.  

\begin{eg}
	Consider a weak Hausdorff space $X$.
	Continuing Example \ref{eg:connect.circulations}, the circulation $\sim$ on $X$ is a k-circulation.  
\end{eg}

\begin{eg}\label{eg:k.d-spaces}
	Continuing Example \ref{eg:d-spaces}, $\leqslant^{dX}$ is a k-circulation for every weak Hausdorff d-space $(X,dX)$.  
\end{eg}

We write $\STREAMS'$ for the category of compactly flowing streams and stream maps between them.
A topological criterion identifies typical examples of such streams.

\begin{prop}\label{prop:topological.criterion}
	Locally compact Hausdorff streams are compactly flowing.
\end{prop}
\begin{proof}
	Consider a locally compact Hausdorff stream $(X,\leqslant)$.
	For each open subspace $U\subset X$, (\ref{eqn:k.cosheaf}) follows from Lemma \ref{lem:pseudo.circulations} because $\mathcal{K}(U)$ contains a neighborhood for each point in its union.
\end{proof}

Our old forgetful functor restricts and corestricts to a new forgetful functor $T':\STREAMS'\rightarrow\SPACES'$ creating final structures exactly as before.

\begin{lem}\label{lem:k.final}
	For each diagram $D:\DIAGRAM\rightarrow\STREAMS'$ and each cocone
	$$\lambda=(\lambda_x:T'Dx\rightarrow X)_{x\in\mathrm{ob}\;\DIAGRAM}$$
	on $T'D$, there is a final structure for $(T',D,\lambda)$ to $(X,\bigvee_{x\in\mathrm{ob}\;\DIAGRAM}\leqslant^{(\lambda_x)_*})$.  
\end{lem}
\begin{proof}
	Trivial circulations on weak Hausdorff spaces are k-circulations.
	The operation $\bigvee$ applied to a non-empty set of k-circulations on a fixed weak Hausdorff space is a k-circulation because $\bigvee$ is commutative.  
	Therefore the result follows from Lemma \ref{lem:final}.
\end{proof}

Consequently, we can characterize compactly flowing streams as quotients of locally compact Hausdorff streams by equivalence relations having closed graphs.

\begin{lem}\label{lem:cogenerators}
	Every compactly flowing stream is the $\STREAMS'$-colimit of compact Hausdorff substreams.
\end{lem}
\begin{proof}
	Consider a compactly flowing stream $(X,\leqslant)$.  
	Let $\DIAGRAM$ denote the category of all compact Hausdorff substreams of $(X,\leqslant)$ and all stream maps between them which are inclusions of underlying sets and let $D:\DIAGRAM\rightarrow\STREAMS'$ denote inclusion.  
	For each $x\in\mathrm{ob}\;\DIAGRAM$, let $\leqslant^x$ denote the circulation of $Dx$.  
	
	The colimit $\bar{\lambda}=\colim\;T'D$ is the cocone on $T'D$ defined by all possible inclusions into $X$, because the diagram $T'D$ contains (at least) all possible inclusions of singletons into compact Hausdorff subspaces of $X$.  
	For each open subspace $U\subset X$,
	\begin{eqnarray}
		\leqslant_U
		&=& \bigvee_{L\in\mathcal{K}(U)}\leqslant_{\restriction L}\label{eqn:cogenerators.k.defn}\\
		&=& \bigvee_{x\in\mathrm{ob}\;\DIAGRAM}\bigvee_{L\in\mathcal{K}(Dx\cap U)}\leqslant_{\restriction L}\label{eqn:cogenerators.k.parse}\\
		&=& \bigvee_{x\in\mathrm{ob}\;\DIAGRAM}\bigvee_{L\in\mathcal{K}(Dx\cap U)}(\leqslant^x)_{\restriction L}\label{eqn:cogenerators.substreams}\\
		&=& \bigvee_{x\in\mathrm{ob}\;\DIAGRAM}\leqslant^x_{Dx\cap U}\label{eqn:cogenerators.compact}\\
		&=& \bigvee_{x\in\mathrm{ob}\;\DIAGRAM}\leqslant^x_{Dx\cap U}\vee\;\mathrm{id}_U\nonumber\\
		&=& \bigvee_{x\in\mathrm{ob}\;\DIAGRAM}(\leqslant^x)^{(\lambda_x)_*}_U\label{eqn:cogenerators.multiple.k.defn},
	\end{eqnarray}
	as (\ref{eqn:cogenerators.k.defn}) and (\ref{eqn:cogenerators.multiple.k.defn}) follow from definitions, (\ref{eqn:cogenerators.k.parse}) follows from $\mathcal{K}(U)=\bigcup_{x\in\mathrm{ob}\;\DIAGRAM}\mathcal{K}(Dx\cap U)$, (\ref{eqn:cogenerators.substreams}) follows from Lemma \ref{lem:inclusion.composites}, and (\ref{eqn:cogenerators.compact}) follows from Proposition \ref{prop:topological.criterion}.
	The final structure for $(T',D,\lambda)$ is a cocone to $(X,\leqslant)$ by Lemma \ref{lem:k.final}.
	Thus $(X,\leqslant)=\colim\;D$.
\end{proof}

\begin{eg}[CW streams]\label{eg:CW}
	For each integer $n\geq 0$, let $\leqslant_{\mathbb{I}^n}$ denote the $n$-fold product of the natural order $\leqslant_{\mathbb{I}}$ on $\mathbb{I}$.
	We can construct stream-theoretic analogues of CW complexes as inductive pushouts of copies of $(\mathbb{I}^n,\leqslant_{\mathbb{I}^n})$, regarded as compactly flowing streams by Theorem \ref{thm:embedding} and Proposition \ref{prop:topological.criterion}, along their boundaries $\partial\mathbb{I}^n$.  
	These ``CW streams'' are compactly flowing by Lemma \ref{lem:cogenerators} and they have CW complexes as their underlying spaces by Lemma \ref{lem:k.final}.  
	Examples include the ``directed circle'' in Example \ref{eg:quotients} and the stream $(\mathbb{R}^2,\leqslant^{d\mathbb{R}^2})$ from Example \ref{eg:d.beaucoup}.  
\end{eg}

As before, a direct consequence of \cite[Propositions 7.3.6 and 7.3.11]{borceux:2} and Lemma \ref{lem:k.final} is that we can also form arbitrary initial structures in $\STREAMS'$.

\begin{prop}\label{prop:k.topological}
	The functor $T':\STREAMS'\rightarrow\SPACES'$ is topological.
\end{prop}

In particular, $T'$ has full and faithful left and right adjoints by \cite[Proposition 7.3.7]{borceux:2} and $T'$ creates limits and colimits by \cite[Proposition 7.3.8]{borceux:2}. 
We must generally ``k-ify'' the cosheafifications of pullbacks in order to form initial structures.
We highlight some cases where k-ification is unnecessary.

\begin{eg}[Compactly generated open substream]
	An open substream of a compactly flowing stream is compactly flowing if its underlying space is compactly generated because the restriction of a k-circulation is a k-circulation.
\end{eg}

\begin{eg}[Closed and ``convex'' substream]
	Consider a compactly flowing stream $(X,\leqslant)$ and a closed substream $(A,\leqslant')$ of $(X,\leqslant)$ such that $A$ is a convex subset of $(U,\leqslant_U)$, for some open neighborhood $U$ of $A$ in $X$.
	The space $A$ is compactly generated.  
	It is straightforward to show $\leqslant'$ is a k-circulation, following our discussion in Example \ref{eg:convex} and Lemma \ref{lem:inclusion.composites}.
\end{eg}

The circulations of (compactly flowing) streams ``preserve finite products.'' 

\begin{prop}\label{prop:products}
	For a $\STREAMS'$-product $(X,\leqslant)=(Y,\leqslant')\times(Z,\leqslant'')$,
	\begin{equation}\label{eqn:products}
	\leqslant_{U\times V}\;=\;\leqslant'_U\times\leqslant''_V
	\end{equation}
	for all open subsets $U\subset Y, V\subset Z$.  
\end{prop}
\begin{proof}
	The space $X$ is the $\SPACES'$-product $Y\times Z$ and projections define stream maps $(X,\leqslant)\rightarrow(Y,\leqslant')$ and $(X,\leqslant)\rightarrow(Z,\leqslant'')$ because $T':\STREAMS'\rightarrow\SPACES'$ preserves products by Proposition \ref{prop:k.topological}.
	Therefore in (\ref{eqn:products}), the graph of the right side contains the graph of the left side; it suffices to show the reverse containment.

	Consider some $y,y'\in Y$ and $z,z'\in Z$ satisfying 
	$$(y,y')\;(\leqslant'_U\times\leqslant''_V)\;(z,z')$$
	for some open subsets $U\subset Y,\;V\subset Z$.
	The identity on $(Y,\leqslant')$ and the constant stream map $(Y,\leqslant')\rightarrow(Z,\leqslant'')$ at $z$ together induce a universal stream map 
	$$\iota:(Y,\leqslant')\rightarrow(X,\leqslant).$$
	The identity on $Y$ and the constant function $Y\rightarrow Z$ at $z$ induce $T'(\iota)$ because $T'$ preserves products. 
	Thus $\iota(y)=(y,z)$ and $\iota(y')=(y',z)$, and hence $(y,z)=\iota(y)\leqslant_{U\times V}\iota(y')=(y',z)$.
	Similarly, $(y',z)\leqslant_{U\times V}(y',z')$.  
	Then $(y,z)\leqslant_{U\times V}(y',z')$ by transitivity.  
	We have thus shown that in (\ref{eqn:products}), the graph of the left side contains the graph of the right side.
\end{proof}


We have constructed a full, complete, and cocomplete subcategory $\STREAMS'\subset\STREAMS$ large enough to contain all streams of interest, including ``mapping streams.''

\begin{thm}\label{thm:xclosed}
	The category $\STREAMS'$ is Cartesian closed.
\end{thm}
\begin{proof}
	Consider a compactly flowing stream $(X_0,\leqslant^0)$.
	The category $\STREAMS'$ has all colimits and finite products by Proposition \ref{prop:k.topological} because $\SPACES'$ has all colimits and finite products.  
	In order to show that the functor
	\begin{equation}\label{eqn:xclosed}
		(X_0,\leqslant^0)\times -:\STREAMS'\rightarrow\STREAMS'
	\end{equation}
	has a right adjoint, it suffices to show that (\ref{eqn:xclosed}) is cocontinuous by the Adjoint Functor Theorem - the solution set condition is trivial because $\SPACES'$ is Cartesian closed and $T'$ is both topological and fibre-small.  
	
	It suffices to take the case $X_0$ compact Hausdorff by Lemma \ref{lem:cogenerators}.
	For brevity, we write $\leqslant^0\!\times\!\leqslant'$ for the circulation of each product stream of the form $(X_0,\leqslant^0)\times(W,\leqslant')$.
	Consider a compactly generated space $Z$.
	The product $X_0\times Z$ has as a basis $\BASIS$ all open subsets of the form $U\times V$ because $X_0$ is compact Hausdorff.  

	Consider a compactly flowing stream $(Y,\leqslant)$ and a continuous function $p:Y\rightarrow Z$.
	The circulations $\leqslant^0\!\times\!\leqslant^{p_*}$ and $(\leqslant^0\!\times\!\leqslant)^{(\mathrm{id}_X\times p)_*}$ coincide because they coincide on $\BASIS$ by Proposition \ref{prop:products} and Lemma \ref{lem:transitive.closure}.  
	Consider a set $\{\leqslant^i\}_{i\in\mathcal{I}}$ of k-circulations on $Z$.  
	The circulations of $\leqslant^0\!\times\!\bigvee_{i\in\mathcal{I}}\leqslant^i$ and $\bigvee_{i\in\mathcal{I}}(\leqslant^0\!\times\!\leqslant^i)$ coincide because they coincide on $\BASIS$ by Proposition \ref{prop:products} and Lemma \ref{lem:demorgan}.
	The functor (\ref{eqn:xclosed}) preserves final structures by Lemma \ref{lem:k.final} and (\ref{eqn:xclosed}) hence preserves colimits because $X_0\times -:\SPACES'\rightarrow\SPACES'$ preserves colimits.
\end{proof}

\section{Future directions for streams}\label{sec:future}
We can systematically enrich spaces from nature - machine state spaces, classifying spaces of small categories, underlying spaces of time-oriented Lorentzian manifolds, one-point compactifications of (ordered) topological vector spaces - with the structure of time and adapt the machinery of homotopy theory to mine such streams for tractable order-theoretic information.  
In our new world, \textit{fundamental categories} from \cite{fgr:ditop}, \textit{homotopy monoids} from \cite{grandis:constructs}, \textit{homology monoids} from \cite{fahrenberg:dihomology, patchkoria:presimplicial}, and \textit{preordered homology groups} from \cite{grandis:ineq} replace fundamental groupoids, homotopy groups, and homology groups.
Potential applications of these tools extend far beyond the static analysis of programs, as suggested in \cite{goubault:guide, gg:s4}.

One possible line of research is a homotopic approach to monoid theory.
Topological simplices naturally admit partial orders, allowing us to lift geometric realization to a functor from simplicial spaces $s\SPACES'$ to $\STREAMS'$ by Theorem \ref{thm:embedding}.  
Our new realization $|B_*(L)|$ of the nerve of a discrete, unital semilattice $L$ has fundamental monoid $L$ (at a natural basepoint) - even though the classifying space $BL$ is contractible.
We believe homotopy invariants on ``classifying streams'' of topological monoids can extract algebraic information inaccessible to classical methods.
Perhaps such invariants can sharpen the necessary homotopic conditions of \cite{brown:geometry, kobayashi:h3, kos:h3, lafont:h3, squier:h3} for monoids to admit finite and complete presentations.  

\section{Acknowledgements}
The author is greatly indebted to the anonymous referee for detailed corrections and suggestions.
Additionally, the author thanks Eric Goubault, Marco Grandis and Peter May for spotting errors and suggesting improvements in earlier drafts.


\begin{thebibliography}{99}
	\bibitem{borceux:2} 
		\newblock {F. Borceux}, 
		\newblock {\em Handbook of Categorical Algebra 2: Categories and Structures}, 
		\newblock {Encyclopedia of Mathematics and its Applications, vol. 51, Cambridge University Press, 1994, pp. xviii+443.}
	\bibitem{brown:geometry}
		\newblock {K. Brown},
		\newblock {\em The Geometry of Rewriting Systems: a proof of the Anick-Groves-Squier theorem},
		\newblock {Algorithms and classification in combinatorial group theory, (Berekely, CA, 1989)},
		\newblock {Math. Sci. Res. Inst. Publ., vol. 23, Springer, New York, 1992, pp. 137-163.}
	\bibitem{bw:models} 
		\newblock {P. Bubenik, K. Worytkiewicz}, 
		\newblock {\em A model category for local pospaces}, 
		\newblock {Homology Homotopy Appl., vol. 8(1), 2006, pp. 263-292.}
	\bibitem{fahrenberg:dihomology} 
		\newblock {U. Fahrenberg}, 
		\newblock {\em Directed Homology},
		\newblock {In Proc. GETCO\&CMCIM 2003, Electronic Notes in Theoretical Computer Science, vol. 100, Elsevier, 2004.}  
	\bibitem{fgrh:components} 
		\newblock {L. Fajstrup, E. Goubault, E. Haucourt, M. Raussen}, 
		\newblock {\em Components of the fundamental category}, 
		\newblock {Appl. Categ. Structures, vol. 12, no. 1, 2004, pp. 84-108.} 
	\bibitem{fgr:ditop} 
		\newblock {L. Fajstrup, E. Goubault, M. Raussen}, 
		\newblock {\em Algebraic topology and concurrency}, 
		\newblock {Theoret. Comput. Sci, vol. 357(1-3), 2006, pp. 241-278.}
	\bibitem{gaucher:flows} 
		\newblock {P. Gaucher}, 
		\newblock {\em A model category for the homotopy theory of concurrency}, 
		\newblock {Homology Homotopy Appl., vol. 5(1), 2003, pp. 549-599.}
	\bibitem{scott:lattices}
		\newblock {G. Gierz, K.H. Hoffman, K. Keimel, J.D. Lawsonj, M. Mislove, and D.S. Scott},
		\newblock {\em Continuous lattices and domains},
		\newblock {vol. 63 of {\em Encyclopedia of Mathematics and Applications}},
		\newblock {Cambridge University Press, Cambridge, 2003.}
	\bibitem{goubault:guide}
		\newblock {E. Goubault},
		\newblock {\em Geometry and Concurrency: A User's Guide},
		\newblock {Mathematical Structures in Computer Sciences, vol. 10, no 4, Aug. 2000, pp. 411-425.}
	\bibitem{gg:s4}
		\newblock {E. Goubault, J. Goubault-Larrecq},
		\newblock {\em On the Geometry of Intutionistic S4 Proofs},
		\newblock {Homology Homotopy Appl., vol. 5(2), 2003, pp. 137-209.}
	\bibitem{gh:components2} 
		\newblock {E. Goubault, E. Haucourt}, 
		\newblock {\em Components of the fundamental category II}, 
		\newblock {Appl. Categ. Structures, vol. 15, no. 4, 2007, pp. 387-414.}
	\bibitem{grandis:d} 
		\newblock {M. Grandis}, 
		\newblock {\em Directed homotopy theory.  I.}, 
		\newblock {Cah. Topol. G\'eom. Diff\'er. Cat\'eg, vol. 44(4), 2003, pp. 281-316.}
	\bibitem{grandis:constructs} 
		\newblock {M. Grandis}, 
		\newblock {\em Directed homotopy theory.  II. Homotopy constructs}, 
		\newblock {Theory and Applications of Categories, vol. 10, no. 14, 2002, pp. 369-391.} 
	\bibitem{grandis:ineq} 
		\newblock {M. Grandis}, 
		\newblock {\em Inequilogical spaces, directed homology and noncommutative geometry}, 
		\newblock {Homology Homotopy Appl., vol. 6(1), 2004, pp. 413-437.}
	\bibitem{grandis:images} 
		\newblock {M. Grandis}, 
		\newblock {\em Ordinary and directed combinatorial homotopy, applied to image analysis and concurrency}, 
		\newblock {Homology Homotopy Appl.}, 
		\newblock {vol. 5(2), 2003, pp. 211-231.}
	\bibitem{haucourt:dicategories} 
		\newblock {E. Haucourt},
		\newblock {\em Comparing topological models for concurrency},
		\newblock {to appear in GETCO 2005 proceedings.}
	\bibitem{kobayashi:h3} 
		\newblock {Y. Kobayashi},
		\newblock {\em Complete rewriting systems and homology of monoid algebras},
		\newblock {Journal of Pure and Applied Algebra, vol. 65 (3), 1990, pp. 263-275.}
 	\bibitem{kos:h3}
		\newblock {Y. Kobayashi, F. Otto, C. Squier},
		\newblock {\em A finiteness condition for rewriting systems},
		\newblock {Theoretical Computer Sciences, vol. 131, 1994, pp. 271-294.}	
	\bibitem{krishnan:thesis}
		\newblock {S. Krishnan},
		\newblock {\em A homotopy theory of locally preordered spaces},
		\newblock {Ph.D. thesis, University of Chicago, Chicago IL, 2006.}  
	\bibitem{lafont:h3}
		\newblock {Y. Lafont}, 
		\newblock {\em A new finiteness condition for monoids presented by complete rewriting systems (after Craig C. Squier)}, 
		\newblock {Journal of Pure and Applied Algebra, Elsevier, vol. 98, 1995, pp. 229-244.}
	\bibitem{nachbin:order} 
		\newblock {L. Nachbin}, 
		\newblock {{\em Topology and order}, translated from the Portugese by Lulu Bechtolsheim},  
		\newblock {Van Nostrand Mathematical Studies, no. 4}, 
		\newblock {D. Van Nostrand Co., Inc., Princeton, N.J.-Toronto, Ont.-London, 1965, vi+122.}
	\bibitem{patchkoria:presimplicial}
		\newblock {A. Patchkoria},
		\newblock {\em Homology and cohomology monoids of presimplicial semimodules},
		\newblock {Bull. Georgian Acad. Sci. vol. 162(1), 2000, pp. 9-12.}
	\bibitem{pratt:models} 
		\newblock {V. Pratt}, 
		\newblock {\em Modelling concurrency with geometry}, 
		\newblock {Proc. 18th ACM Symp. on Principles of Programming Languages, ACM Press, New York, 1991, pp. 311-322.} 
	\bibitem{squier:h3}
		\newblock {C. C. Squier},
		\newblock {\em Word problems and a homological finiteness condition for monoids},
		\newblock {Journal of Pure and Applied Algebra, vol. 49, 1987, pp. 201-217.}
	\bibitem{steenrod:convenience} 
		\newblock {N. Steenrod}, 
		\newblock {\em A convenient category of spaces},
		\newblock {Michigan Math. J., vol. 14, 1967, pp 133.152.}
\end{thebibliography}
\end{document}